\documentclass [12pt] {article}
\usepackage{amsfonts}
\usepackage{amssymb}

\newcommand{\qed} {\hspace {0.1in} \rule {1.5mm} {3.5mm}}

\newtheorem{lemma}{Lemma}[section]
\newtheorem{corollary}{Corollary}[section]
\newtheorem{theorem}{Theorem}

\newtheorem{proposition}{Proposition}[section]
\newtheorem{definition}{Definition}[section]
\def\a{\alpha}

\def\limn{\lim_{n\to\infty}}
\def\limk{\lim_{k\to\infty}}

\def\e{\epsilon}

\def\deg{\mbox{deg}\,}

\def\tr{\mbox{Tr}\,}
\def\supp{\mbox{supp}\,}
\def\mat{\mbox{Mat}\,}
\def\proj{\mbox{proj}\,}
\def\matd{\mat_{d\times d}}

\def\dim{{\rm dim}}
\def\dimc{\dim_{\bC}}
\def\supp{{\rm supp}}

\def\<{\langle}
\def\>{\rangle}

\def\proof{\smallskip\noindent{\it Proof.} }

\def\bZ{{\mathbb Z}}
\def\bR{{\mathbb R}}
\def\bN{{\mathbb N}}
\def\bC{{\mathbb C}}
\def\deg{\mbox{deg}\,}

\def\to{\rightarrow}

\begin{document}
\title{The Strong Approximation Conjecture holds for amenable groups} 
\author{\sc G\'abor Elek
\footnote {The Alfred Renyi Mathematical Institute of
the Hungarian Academy of Sciences, P.O. Box 127, H-1364 Budapest, Hungary.
email:elek@renyi.hu}}
\date{}
\maketitle
\vskip 0.2in
\noindent{\bf Abstract.}\hskip 0.2in
Let $G$ be a finitely generated group and $G\rhd G_1\rhd
G_2\rhd \dots $ be normal subgroups such that
 $\cap_{k=1}^\infty G_k=\{1\}$. Let
$A\in \matd(\bC G)$ and $A_k\in \matd(\bC ( G/G_k))$ be the images of $A$
under the maps induced by the epimorphisms $G\to G/G_k$. According to the
strong form of the Approximation Conjecture of L\"uck \cite{Luck}
$$\dim_G(\ker\, A)=\lim_{k\to\infty} \dim_{G/G_k}(\ker\, A_k)\,,$$ where
$\dim_G$ denotes the von Neumann dimension.
 In \cite{DLM} Dodziuk, Linnell, Mathai,
Schick and Yates proved the conjecture for torsion free elementary amenable
groups. In this paper we extend their result for all amenable groups, using
the quasi-tilings of Ornstein and Weiss \cite{OW}.
\vskip 0.3in
\newpage
\section{Introduction}

First, let us recall the approximation result of Dodziuk, Linnell, Mathai, 
Schick and Yates \cite{DLM}. Let $G$ be a finitely generated
group and let $A\in \matd(\bC G)$ .

\noindent
Let $l^2(G)=\{f: G\to \bC\,\mid \sum_{g\in G} 
|f(g)|^2<\infty\}$. By left convolution, $A$ induces a
bounded linear operator $A:(l^2(G))^d\to  (l^2(G))^d$,
which commutes with the right $G$-action. Let
$$\proj_{\ker\,A}:(l^2(G))^d\to  (l^2(G))^d$$
be the orthogonal projection onto $\ker\,A$. Then
$$\dim_G(\ker\,A):=\tr_G(\proj_{\ker\,A}):=
\sum^d_{i=1}\langle \proj_{\ker\,A}\,{\bf 1}_i,{\bf 1}_i
\rangle_{(l^2(G))^d}\,$$
where ${\bf 1}_i\in (l^2(G))^d$ is the function which takes the value
 $e_i$ on the 
unit element of $G$  and zero elsewhere ($\{e_1, e_2,\dots, e_n\}$ is an 
orthonormal basis of $\bC^d$).
$\dim_G (\ker\,A)$ is called the von Neumann dimension of
$\ker\,A$.

\noindent
Now let $G\rhd G_1\rhd
G_2\rhd\dots $ be normal subgroups such that
 $\cap_{k=1}^\infty G_k=\{1\}$.

\noindent
Let $A_k\in \matd(\bC ( G/G_k))$ be
 the images of $A$
under the maps induced by the epimorphisms $G\to G/G_k$. According to the
strong form of the Approximation Conjecture of L\"uck \cite{Luck}
$$\dim_G(\ker\, A)=\lim_{k\to\infty} \dim_{G/G_k}(\ker\, A_k)\,.$$ 
In \cite{DLM} the authors prove the conjecture above in the case
when $G$ is a torsion-free elementary amenable group. The goal
of this paper is to extend their result to arbitrary amenable groups.
If $A\in \matd(\bZ(G))$ the problem is much easier to handle 
since one can use the method of L\"uck \cite{Luck2}. Then the conjecture
holds for a large class of groups including amenable and residually
finite groups. In the case of complex group algebra the situation seems
much more complicated. Dodziuk et al. used noncommutative algebra to prove the
conjecture, we
shall use the quasi-tilings of Ornstein and Weiss.

\section{Preliminaries}
Let $G$ be a finitely generated amenable group with a finite symmetric
set of generators $S$. Consider the Cayley-graph $C_G$, where $V(C_G)=G$
and 
$$E(C_G):=\{(x,y)\in G\times G\,\mid\, y=sx, s\in S\}\,.$$
Now we introduce some notation frequently used in the paper later on.
\begin{enumerate}
\item
If $g\in G$, then its word-length $w(g)$ is defined as $d_{C_G}(g,1)$,
where $d_{C_G}$ is the shortest path distance on the Cayley-graph.
\item
Let $F\subset G$ be a finite set, $k>0$, then $B_k(F)$ denotes
the $k$-neighborhood of $F$ in the $d_{C_G}$-metric.
\item
We denote by $\Omega_k(F)$ the set of vertices $p$ in $F$, such that
$d_{C_G}(p,F^c)>k$, where $F^c$ is the complement of $F$.
\item
For $A\in \matd(\bC G)$, its propagation $w(A)$ is just $\sup\,w(g)$,
where $g$ runs through the terms of
non-zero coefficients in the entries of $A$.
Observe that if $f\in (l^2(G))^d$ and $\supp f=U\subset G$, then
$\supp\, A(f)\in B_{w(A)}(U)$, and if $\supp f=\Omega_{w(A)} (U)$ then
$\supp\,A(f)\subseteq U$. 
\item
For a finite set $F\subset G$, $\partial F$ denotes the
set of vertices in $F$ such that $d_{C_G} (p,F^c)=1$. We shall denote the
ratio $\frac {|\partial F|} {|F|}$ by $i(F)$.
\item
Since $G$ is amenable, it has a {\it Folner-exhaustion}, that is
a sequence of subsets
$1\in F_1\subset F_2\subset\dots, \cup_{n=1}^\infty F_n=G$ such that
$i(F_n)\to 0$.
\end{enumerate}

 Now we prove some approximation theorems for amenable groups. 
 Let $1\in F_1\subseteq F_2\subseteq\dots $ be a
Folner exhaustion of $G$ and $P_n:(l^2(G))^d\to (l^2(F_n))^d$ be the 
orthogonal projections.
Then by 
{\it Theorem $3.11$} \cite{DLM} (or Proposition 1. \cite{Elek}) :
$$\dim_G(\ker\, A)=\lim_{n\to\infty}
\frac{\dim_{\bC} (\ker\,P_nAP_n^*)}{|F_n|}\,.$$
\noindent
We define the following sequences of vector spaces:
$$Z_n:=\{f\in (l^2(G))^d\,\mid\, \supp f\subseteq 
B_{w(A)}(F_n), A(f)\mid_{F_n}=0\}$$
$$W_n:=\{f\in (l^2(G))^d\,\mid\, \supp f\subseteq \Omega_{w(A)}(F_n)
  , A(f)=0\}$$
$$V_n:=\ker\,(P_nAP_n^*)$$
\begin{proposition} \label{p5} $$
\limn \frac{\dimc (Z_n)}{|F_n|}=\dim_G (\ker A),\quad
\limn \frac{\dimc (W_n)}{|F_n|}=\dim_G (\ker A). $$
\end{proposition}
\proof
It is enough to prove that
\begin{equation} \label{e5}
\limn \frac{\dimc (Z_n)}{\dimc (V_n)}=1
\end{equation}
and
\begin{equation} \label{e51}
\limn \frac{\dimc (W_n)}{\dimc (V_n)}=1
\end{equation}
Clearly, $W_n=V_n\cap \{f\in l^2(F_n)^d\,\mid\, 
\supp f\subseteq \Omega_{w(A)}(F_n)  \}.$ 
Hence (\ref{e51}) follows from the fact that
$$\limn \frac {|\Omega_{w(A)}(F_n)|} {|F_n|}=1\,.$$
Also, $W_n=Z_n\cap \{f\in l^2(F_n)^d\,\mid\, 
\supp f\subseteq \Omega_{w(A)}(F_n)  \}.$ 
Therefore (\ref{e5}) follows from the fact that
$$\limn \frac {|\Omega_{w(A)}(F_n)|} {|B_{w(A)}(F_n)| }=1\,\quad \qed$$

\noindent
\begin{definition} \label{prop}
Let $F\subset G$ be a finite set and $\delta, \epsilon >0$ be real
numbers. We say that $F$ has property $A(\epsilon, \delta,-)$
if for {\it any} subset $K\subseteq F$, $\frac{|K|} {|F|}> 1-\epsilon$,
the following holds:
\begin{itemize}
\item If $R=\{f\in l^2(F_n)^d\,\mid\, 
\supp f\subseteq \Omega_{w(A)}(K), A(f)=0\}$, then

\noindent
$\dimc\,(R)\geq (1-\delta)
(\dim_G(\ker A))\,.$
\end{itemize}
Also, we say that $F$ has property $A(\delta,+)$
if the following holds:
\begin{itemize}
\item
If $Q$ is the restriction of the space
$$Z_F:=\{f\in (l^2(G))^d\,\mid\, \supp f\subseteq 
B_{w(A)}(F), A(f)\mid_{F}=0\}$$
onto $F$, then
$$\dimc (Q)\leq (\dim_G(\ker A)) +\delta$$. \end{itemize}
\end{definition}
Similarly to Propostion \ref{p5} one can easily prove the following
proposition.
\begin{proposition} \label {p5b} Let 
$1\in F_1\subseteq F_2\subseteq\dots $ be a Folner exhaustion of $G$
as above. Then
for any pair of real numbers $\delta, \epsilon >0$ 
there exists $n_{\delta, \epsilon}$ such that if
$n\geq n_{\delta, \epsilon}$ then 
$F_n$ has both properties $A(\epsilon, \delta,-)$ and
$A(\delta,+)$.
\end{proposition}

\section{Graph convergence and dimension averaging} \label{s2}

Let $C_G$ be the Cayley-graph of the previous section. Color the
directed edge $(x\to y), x=sy$ by $s\in S$ (hence $(y\to x)$ shall
be colored by $s^{-1}$). Thus we color all edges in both direction
with the elements of the set $S$ such a way that for each
$x\in G$ the edges outgoing from $x$ are colored in different ways. 
The following definition is a variation of the one on random weak
convergence in \cite{BS}.

\noindent
Let $B_1,B_2,\dots$ be an infinite sequence of finite graphs. Assume
that for any $x\in V(B_n)$: $\deg(x)\leq |S|$.
We also assume that the directed edges
are colored by $S$ such a way that :
\begin{itemize}
\item the color of the edge $(x\to y)$ is the inverse of the color of 
$(y\to x)$.
\item the outgoing edges from any vertex are colored differently.
\end{itemize}
We say that $p\in V(B_n)$ is $k$-similar to the identity of $G$,
if its $k$-neighborhood in $B_n$ is edge-colored isomorphic to
the $k$-neighborhood of the identity in $C_G$. Let $Q^B_{k}$ be the 
set of vertices in $B$ that are $k$-similar to the identity. Then we say that
$\{B_n\}^\infty_{n=1}$ converge to $C_G$ if for any $\epsilon>0$ and
$k\in \bN$ there exists $n_{\epsilon,k}$ such that if $n\geq n_{\epsilon,k}$
then 
$$ Q^{B_n}_{k}>(1-\e) |V(B_n)|\,.$$

\noindent
{\bf Example 1.:\,} Let $G$ be a finitely generated group
and $\{B_n\}^\infty_{n=1}$ be a sequence of finite induced subgraphs
forming a Folner-exhaustion. Then $\{B_n\}^\infty_{n=1}$ converge to
$C_G$.

\noindent
{\bf Example 2.:\,} Let $G$ be a finitely generated residually
finite group and $G\rhd G_1 \rhd G_2 \rhd\dots...$ be
a sequence of finite index normal subgroups such that
$\cap^\infty_{n=1} G_n=\{1\}$. Let $C_n$ be the Cayley-
graph of $G/{G_n}$. Then $\{C_n\}^\infty_{n=1}$ converge to
$C_G$.

\noindent
Now let $A\in \matd(\bC G)$. One can define the transformation kernel
of $A$, $\widetilde{A}:G\times G\to \matd(\bC)$ the following way.
First write $A$ in the form of $\sum_{\gamma\in G} A_{\gamma}\cdot \gamma$,
where $A_\gamma\in \matd (\bC)$. Then set $\widetilde{A}(\gamma, \delta):=
A_{\gamma\delta^{-1}}\,.$ Thus if $f\in l^2(G))^d$,
then
$$A(f) (\delta)=\sum_{\gamma\in G} \widetilde{A}(\delta, \gamma) f(\gamma)\,.$$
Now let 
$\{B_n\}^\infty_{n=1}$
be a sequence of graphs converging to $C_G$.
Then we define the
finite dimensional linear
transformations $T_n:(l^2(V(B_n))^d\to (l^2(V(B_n))^d$
approximating $A$, the following way.
\begin{itemize}
\item
If $ x\in Q^{B_n}_{w(A)}$, $y\in V(B_n)$ and
$d_{B_n}(y,x)\leq w(A)$,
let
$\widetilde{T}_n(y,x):=A(\gamma,1)$, where
 $\gamma$ is the element of
$G$ satisfying $\phi^x_{w(A)}(\gamma)=y$. Here $\phi^x_{w(A)}$ is the unique
 colored isomorphism between the $w(A)$-neighborhood of $1$ in $C_G$ and
the $w(A)$-neighborhood of $1$ in $B_n$.
\item
If $x\notin Q^{B_n}_{w(A)}$ or $d_{B_n}(y,x)>w(A)$, then
let $\widetilde{T}_n(x,y):=0$. \end{itemize}
Then if $f\in l^2(V(B_n))^d$ and $p\in V(B_n)$;
$$T_n(f) (p)=\sum_{q\in V(B_n)} \widetilde{T}_n(p,q) f(q)\,.$$

\noindent
The main goal of our paper is to prove the following theorem.
\begin{theorem} \label{t11}
If $G$ is a finitely generated amenable group and
$\{B_n\}_{n=1}^\infty$, $\{T_n\}_{n=1}^\infty$ are as above, then

$$\limn \frac {\dim_{\bC} \ker T_n}{|V(B_n)|}=\dim_G(\ker A)$$
\end{theorem}
The Strong Approximation Conjecture for amenable groups follows
from the theorem:

\begin{corollary}
If $G$ is a finitely generated amenable group and \\
$G\rhd G_1 \rhd G_2\dots,\cap^\infty_{n=1} G_n=\{1\}$ are
normal subgroups, then
$$\limn \dim_{G/G_n}(\ker A_n)= \dim_G(\ker A)\,,$$
where $A\in \matd(\bC G)$ and $A_n\in \matd(\bC ( G/G_n))$
are the images of $A$
under the maps induced by the epimorphisms $G\to G/G_n$. \end{corollary}
\proof (of the Corollary)

\noindent
\underline{Case 1.} Suppose that all $G_n$ has finite index. Note that
in this case $T_n=A_n$ if $n$ is large enough, hence the corollary immediately
follows.

\noindent
\underline{Case 2.} Assume that for large enough $n$, the
amenable group $G/G_n$ is infinite.
Let $1\in F^n_1\subset F^n_2\subset\dots$ be a Folner-exhaustion of the Cayley
graph
$C_{G/G_n}$ (using the image of the generator system $S$).
Then
$$\dim_{G/G_n}(A_n)=\limk \frac{\dim_{\bC} (\ker P^n_k A_n (P^n_k)^*)}
{|F^n_{k}|}\,,$$
where $P^n_k:(l^2(G/G_n))^d\to (l^2(F^n_{k}))^d$ is
the orthogonal projection.

\noindent
Pick a sequence $F^1_{m_1}, F^2_{m_2},\dots$ such that
\begin{itemize}
\item $i(F^j_{m_j})\to 0$.
\item $(\dim_{G/G_n}(A_n)- \frac{\dim_{\bC} (\ker P^{n}_{m_n} A_n
  (P^n_{m_n})^*)}{|F^n_{m_n}|})\to 0$.
\end{itemize}
Now let $B^n_{m_n}$ be the graph induced by $F^n_{m_n}$. 
\begin{lemma} \label{l13}
$\{B^n_{m_n}\}^\infty_{n=1}$ converge to $C_G$.
\end{lemma}
\proof
Since $\cap^\infty_{k=1} G_k=\{1\}$, for any $d\in\bN$ there exists $n_d>0$
such that if $n\geq n_d$ then
 the $d$-balls in $C_{G/G_n}$ are colored-isomorphic
to
the $d$-ball of $C_G$.
Let $H^n_{m_n}=\Omega_d (F^n_{m_n})$. Clearly $H^n_{m_n}\subseteq 
Q^{F^n_{m_n}}_{d}$.
  Since the vertex degrees of $G/G_n$ are
at most $S$, $|H^n_{m_n}|\geq |F^n_{m_n}|- |S|^d |\partial F^n_{m_n}|\,.$
Now our lemma easily follows. \qed
\begin{lemma}
\label{l14}
$$\limn \frac {\dim_{\bC} (\ker P^n_k A_n (P^n_k)^*)} {\dim_{\bC} \ker
  T_n}=0\,.$$
Here $T_n$ is the linear operator associated to $B^n_{m_n}$.
\end{lemma}
If $\supp f\subset F^n_{m_n}\backslash B_{w(A)}(\partial F^n_{m_n})$ then
$T_n(f)= P^n_k A_n (P^n_k)^* (f)\,.$
Since
$$\frac {| F^n_{m_n}\backslash B_{w(A)}(\partial F^n_{m_n})|}{|F^n_{m_n}|}\to
1\,$$
our lemma follows. \qed

\noindent
Obviously, Lemma \ref{l13} and Lemma \ref{l14} imply the corollary.
\qed
\section{Quasi-tilings}
Let us recall the notion of quasi-tilings from \cite{OW}.
Let $X$ be a finite set and $\{A_i\}^n_{i=1}$ be finite subsets of $X$. Then
we say that $\{A_1, A_2,\dots,A_n\}$ are $\epsilon$-disjoint if there exist
subsets $\overline{A_i}\subseteq A_i$ such that
\begin{itemize}
\item For any $1\leq i \leq n$, $\frac {|\overline{A_i}|}{|A_i|}\geq
  1-\epsilon.$
\item If $i\neq j$ then $\overline{A_i}\cap \overline{A_j}=\emptyset$.
\end{itemize}
On the other hand, if $\{H_j\}^m_{j=1}$ are finite subsets of $X$, then we say
that
they $\alpha$-cover $X$ if
$$\frac{|X\cap (\cup^m_{j=1} H_j)|} {|X|}\geq \alpha\,.$$
Finally, we say that the collection $\{H_1, H_2,\dots,H_m\}$ $\delta$-evenly
covers $X$
if there exists some $M\in\bR^+$ such that 
\begin{itemize}
\item For any $p\in X$, $\sum^m_{j: p\in H_j} 1\leq M$.
\item $\sum^m_{j=1} |H_j| \geq (1-\delta) M |X|\,.$
\end{itemize}
According to Lemma 4. \cite{OW}, if $\{H_1,H_2,\dots,H_m\}$ form
a $\delta$-even covering of $X$, then for any 
$0<\epsilon <1$ there exists an $\e$-disjoint subcollection of
the $H_j'$s that $\e(1-\delta)$-covers $X$.

\noindent
Now we define {\it tiles} for our $S$-edge colored graphs. Let $G$ be a
finitely generated group with a symmetric generator set $S$ and let
$1\in F_1\subseteq F_2\subseteq\dots\,, \cup^\infty_{n=1} F_n=G$ be
a Folner-exhaustion.
Let $B$ be a finite graph as in the previous section with edge-colorings
by the elements of $S$. Also, let $L$ be a natural number.
Let $\{F_{\alpha_1}, F_{\alpha_2},\dots F_{\alpha_n} \}$ be a finite
collection of the Folner sets above such that for
any $1\leq i \leq n$, $F_{\alpha_i}\subset 
B_{\frac{1}{2}L}(1)$. Then for any $x\in Q^B_L$ and $1\leq i \leq n$, 
$T_x(F_{\alpha_i})$ is the image of $F_{\alpha_i}$ under the unique
colored isomorphism $\phi^x_L:B_L(1)\to B_L(x)$ mapping $1$ to $x$.
We call such a subset a tile of type $F_{\alpha_i}$ and say that
$x$ is the center of $T_x(F_{\alpha_i})$.
A system of tiles $\e$-quasi tile $V(B)$ if they are $\e$-disjoint and form
an $(1-\e)$-cover.
The following theorem is a version of Theorem 6. in \cite{OW}.
\begin{theorem}
\label{t17}
For any $\e>0$, $n>0$, there exist $L>0$, $\delta>0$ and a finite
collection $\{F_{n_1}, F_{n_2},\dots, F_{n_s}\}\subset B_L(1)$
 of the Folner sets, such that $n_i>n$ and if
$$\frac {|Q^B_L|}{|V(B)|}>1-\delta$$ then $V(B)$ can be 
$\e$-quasi-tiled by tiles of the form $T_x(F_{n_i})$,
$x\in Q^B_L$, $1\leq i \leq s$.
\end{theorem}
\section{The inductional step}
First of all fix a constant $\e_1<\frac{\e}{100}$. 
Let us call a finite set $H\subset G$ a set of type $(K,\alpha)$,
$K\in \bN$, $\a\geq 0$ if
\begin{equation} \label {1e18}
\frac{|B_K(H)|} {|H|} <1+\a\,. \end{equation}
Now let $B$ be our $S$-edge colored finite graph and
suppose that
\begin{equation} \label {2e18}
\frac {|Q^B_L|}{|V(B)|}>1-\beta\,.
\end{equation}
The exact values of $\beta$ and $L$ shall be given later. Assume
that $H$ is of type $(K,\a)$, where

\begin{equation} \label {3e18}
H\subset B_{\frac{1}{100}L}(1)\quad\mbox{and}\quad K<\frac{L}{10}\,.
\end{equation}

 Now consider all tiles in $B$ in the form $T_x(H)$, where
 $x\in Q^B_L$. Note that no vertices of $B$ is covered by more than $|H|$
tiles. Indeed, if $z$ is covered, then the $\frac{L}{2}$-neighborhood of
$z$ in $B$ is colored isomorphic to the $\frac{L}{2}$-neighborhood of $1$ in
$G$.
Hence if $z\in T_x(H)$, then $z\in Q^B_{\frac{L}{2}}$ and $x\in T_z(H^{-1})$. 
Summarizing these:
\begin{itemize}
\item For any $y\in V(B)$, $\sum_{x: y\in T_x(H)} 1\leq |H|$.
\item $\sum_{x\in Q^B_L} |T_x(H)|= |Q^B_L||H|\geq (1-\beta) |V(B)| |H|\,.$
\end{itemize}
Consequently, the tiles $\{T_x(H)\}_{x\in Q^B_L}$
form a $\delta$-even covering of $V(B)$, where $\delta=1-\beta$.
Then by Lemma 4. of \cite{OW}, there exists an $\e_1$-disjoint 
subcollection of tiles, $\cup_{x\in I} T_x(H)$ such that they
form a $\e_1 (1-\beta)$-covering of $V(B)$.

\noindent
Now suppose that the number of vertices in $V(B)$ not covered by this
subcollection above is greater than $\frac{\e}{2} |V(B)|$.
Let $B_1$ be the graph induced by the uncovered vertices. We would like
to estimate the quotient:
$\frac {|Q^{B_1}_K|}{|V(B)|}\,.$
Note that if $y\in V(B_1)$ and 
$$y\notin \cup_{x\in I} (T_x(B_K(H))\backslash T_x(H))$$
then $y\in Q^{B_1}_K\,.$
Hence by $\e_1$-disjointness,
$$|\cup_{x\in I} (T_x(B_K(H))\backslash T_x(H))|\leq
\a\sum_{x\in I} |H|\leq \a(1-\e_1)^{-1}|V(B)|\,.$$
Hence
$$ |Q^{B_1}_K|\geq |V(B_1)|-\a (1-\e_1)^{-1} |V(B)|\,,$$
that is
\begin{equation}\label{e21}
\frac{|Q^{B_1}_K|} {|V(B)|}\geq 1-\beta_1\,,\end{equation}
where
$\beta_1=\a (1-\e_1)^{-1} \frac {2} {\e}\,$.
Also note that
$$\frac {\e} {2} |V(B)|\leq |V(B_1)| \leq (1-\e_1(1-\beta)) |V(B)|\,.$$

\section{The proof of Theorem \ref{t17}}
Let $\{\a_k\}^\infty_{k=1}$ be a sequence of real numbers tending to zero
and let $\{s_k\}^\infty_{k=1}$ be a sequence of real numbers tending to 
infinity, satisfying the following inequalities:
$$s_k\geq 1, s_{k+1}\geq 10 s_k\,.$$
We call a subsequence of the Folner exhaustion $\{F_n\}^\infty_{n=1}$
an $(\a, s)$-good subsequence if it satisfies the following conditions:
\begin{itemize}
\item
$1\in F_{n_1} \subset B_{s_1}(1) \subset F_{n_2}
\subset B_{s_2}(1) \subset F_{n_3}\subset\dots$
\item
$F_{n_{i+1}}$ is of type $(100s_i,\a_i)$.
\end{itemize}
Obviously  one can choose $\{s_k\}^\infty_{k=1}$ for any fixed
$\{\a_k\}^\infty_{k=1}$ to have such  $(\a, s)$-good subsequences.
Now let $M$ be an integer such that
\begin{equation}
\label{1e22}
(1-\frac{\e_1}{2})^M <\frac {\e} {100}
\end{equation}
Also, pick $\beta>0$ so that
\begin{equation}
\label{2e22}
\beta M< \frac {\e} {100}
\end{equation}
And finally fix a sequence $\{\a_k\}^\infty_{k=1}$ such that
\begin{equation}
\label{3e22}
\a_i (1-\e_1)^{-1} \frac {2} {\e}<\beta\,.
\end{equation}
Now let $B$ a finite $S$-colored graph such that
$$\frac{Q^B_{100 s_M}}{|V(B)|}> 1-\beta\,,$$
where $\beta, M$ are as above. Then by the argument of the
previous section we can $\e_1(1-\beta)$-cover
the vertices of $B$ by $\e_1$-disjoint tiles of type $F_{n_M}$.
If $B_1$ is the graph induced by the uncovered vertices of $V(B)$,
by (\ref{e21}):
$$\frac{|Q^{B_1}_{s_M}|} {|V(B_1)|} > 1-\beta_1\,,$$
where $\beta_1=\a_M (1-\e_1)^{-1} \frac {2} {\e}\,$.
Now we can $\e_1(1-\beta_1)$-cover $V(B_1)$ by tiles of type $F_{n_{M-1}}$.
If $B_2$ denotes the graph induced by the uncovered part of $V(B_1)$ then
$$\frac{|Q^{B_2}_{s_{M-1}}|} {|V(B_2)|} > 1-\beta_2\,,$$
where $\beta_2=\a_{M-1} (1-\e_1)^{-1} \frac {2} {\e}\,$.

\noindent
We proceed inductively. In each step the new tiles are disjoint
from all previous ones. Also,
$$V(B_i)\leq V(B_{i-1}) (1-\frac{\e_1} {2})\,.$$
Hence by our conditions, in at most $M$ steps we obtain
an $\e$-disjoint $(1-\e)$-covering of $V(B)$. \qed

\section{The proof of Theorem \ref{t11}}
Let $G$ be a finitely generated amenable group, $A\in \matd(\bC G)$
and $\{B_n\}^\infty_{n=1}$ be a sequence converging to $C_G$. Let 
$\{T_n\}^\infty_{n=1}$ be the sequence of approximating operators
as in Section \ref{s2}.
\begin{proposition} \label{p25}
For any pair $\delta, \e>0$ there exists $k_{\delta,\e}>0$ such
that if $k\geq k_{\delta,\e}$ then
$$\frac {\dimc\, (\ker T_k)} {|V(B_k)|}\geq
(\dim_G (\ker A)-\delta) (1-\e)\,.$$ \end{proposition}
\proof
Let $1\in F_1\subseteq F_2\subseteq\dots$ be a Folner exhaustion of $G$,
such that all the $F_n$'s have property $A(\e,\delta,-)$
 (see Lemma \ref{prop}).
Let $\{H_1, H_2,\dots H_s \}$ be an $\e$-quasi-tiling of $B_k$ by tiles
from this Folner sequence. Such $\e$-quasi-tiling exists by Theorem \ref{t17}
if $k$ is large enough.
For $1\leq i \leq s$ let $K_i\subset H_i$ be a subset such that
\begin{itemize}
\item $\frac{|K_i|}{|H_i|}> 1-\e\,.$
\item $K_i\cap K_j=\emptyset$ if $i\neq j$.
\end{itemize}
Since the $F_n$'s have property $A(\e,\delta,-)$ there exist
subspaces $V_i\subset (l^2(B_n))^d$ such that
\begin{itemize}
\item If $f\in V_i$ then $\supp\, f\subseteq K_i$.
\item $f\in \ker T_k$.
\item $\frac{\dimc V_i} {|H_i|} \geq \dim_G (\ker A)-\delta\,.$
\end{itemize}
Now consider the subspace $\oplus^s_{i=1} V_i\subseteq \ker T_k\,.$
Then
$$\dimc (\oplus^s_{i=1} V_i)\geq (\sum^s_{i=1} |H_i|)
 (\dim_G (\ker A)-\delta  )\geq
(1-\e) |V(B_k)| (\dim_G (\ker A)-\delta  )\,.$$
That is
$$\frac {\dimc\, (\ker T_k)} {|V(B_k)|}\geq
(\dim_G (\ker A)-\delta) (1-\e)\,. \quad\qed$$

\begin{proposition} \label{p28}
For any pair $\delta, \e>0$ there exists $m_{\delta,\e}>0$ such
that if $k\geq m_{\delta,\e}$ then
$$\frac {\dimc\, (\ker T_k)} {|V(B_k)|}\leq (1-\e)^{-1}
(\dim_G (\ker A)+\delta)+\e\,.$$ \end{proposition}
\proof
Again let $1\in F_1\subseteq F_2\subseteq\dots$ be a Folner exhaustion of $G$,
such that all the $F_n$'s have property $A(\delta,+)$
 (see Lemma \ref{prop}).
Consider the $\e$-quasi-tilings of the previous proposition.
Now let $W_i\subset l^2(H_i)$ be the restriction
of $\ker T_k$ onto $H_i$.
By our assumption,
$$\dim_G(\ker T_k) \leq |H_i| (\dim_G (\ker A)+\delta)\,.$$
Since $\{H_i\}_{i=1}^s$ form an $\e$-covering 
$$\dimc\, (\ker T_k)\leq \e |V(B_k)| + \sum^s_{i=1} |H_i|
 (\dim_G (\ker A)+\delta)\,.$$
Note that by $\e$-disjointness
$$\sum^s_{i=1} |H_i|\leq (1-\e)^{-1} |V(B_k)|\,.$$
Thus
$$\frac {\dimc\, (\ker T_k)} {|V(B_k)|}\leq (1-\e)^{-1}
(\dim_G (\ker A)+\delta)+\e\,.\quad\qed $$ 
Clearly, Propositions \ref{p25} and \ref{p28} imply Theorem \ref{t11}. \qed

\end{document}